\documentclass[12pt]{amsart}
\usepackage[utf8]{inputenc}
\usepackage[T2A]{fontenc}
\usepackage{textcase}

\usepackage[all, cmtip]{xy}
\usepackage{amsfonts,amsmath,amstext,amsbsy,amssymb, amsopn,amsthm,upref,eucal, mathrsfs}

\newcommand{\dLP}{d_{\mathrm{LP}}}

\newcommand{\ind}{\mathbf{1}}

\newcommand{\PP}{\mathbb{P}}
\newcommand{\EE}{\mathbb{E}}
\newcommand{\NN}{\mathbb{N}}
\newcommand{\ZZ}{\mathbb{Z}}
\newcommand{\RR}{\mathbb{R}}
\newcommand{\T}{\mathbb{T}}
\newcommand{\cH}{\mathcal{H}}
\newcommand{\cM}{\mathfrak{M}}
\newcommand{\cX}{\mathscr{X}}
\newcommand{\rX}{\mathscr{X}}
\newcommand{\dH}{\dot{H}_{-1/2}(\T)}
\newcommand{\ddH}{\dot{H}_{-1/2}(\T)^{\mathbb R}}
\newcommand{\CC}{\mathbb{C}}

\newtheorem{lemma}{Lemma}
\newtheorem{corollary}{Corollary}

\newtheorem{assumption}{Assumption}
\newtheorem{definition}{Definition}

\newtheorem{proposition}{Proposition}
\theoremstyle{remark}

\begin{document}

\title[Convergence in Distribution of Random Fields]{Convergence in Distribution of Random Fields Taking Values in Sobolev Spaces}

\author{Alexander\,I.~Bufetov}
\address{Faculty of Mathematics and Computer Science, Saint Petersburg State University\\
Steklov Mathematical Institute of the Russian Academy of Sciences}
\date{}

\begin{abstract}
The aim of this paper is to develop a general formalism for establishing convergence in distribution of random fields taking values in Sobolev spaces. In particular, the formalism is applied to proving the convergence in distribution of random measures to Gaussian Multiplicative Chaos.
\end{abstract}

\thanks{This work was supported by the Ministry of Science and Higher Education of the Russian Federation under Agreement No.~075-15-2024-631, as part of the implementation of a scientific project.}
\maketitle

\section{Introduction}

The convergence of a sequence of random measures $M^{(n)}$, defined
on a common standard Borel probability space
$(\Omega,\mathcal{B},\PP)$ and taking values in the space of finite
Borel measures on a complete separable metric space $S$, is often
established using the following simple observation. Suppose that a measure $ds$ is given on $S$, and, in addition, that we have
a sequence of Borel kernels $Q^{(n_1,n_2)}$, $n_1,n_2\in\NN$,
defined on $S\times S$ and satisfying the identity
\[
  \EE_{\PP}M^{(n_1)}(I)M^{(n_2)}(J)
  =\iint\limits_{I\times J}Q^{(n_1,n_2)}(s_1,s_2)\,ds_1\,ds_2.
\]
for all Borel subsets $I,J\subset S$.
In this case, convergence of the kernels $Q^{(n_1,n_2)}$ jointly
with respect to the indices $(n_1,n_2)$ implies convergence of the
random measures $M^{(n)}$ as well: for example, if, for a Borel subset
$I\subset S$, the limit
\[
  \lim_{n_1,n_2\to\infty}\iint\limits_{I\times I}Q^{(n_1,n_2)}(s_1,s_2)\,ds_1\,ds_2,
\]
exists, then the sequence $M^{(n)}(I)$ converges in
$L_2(\Omega,\PP)$ as $n\to\infty$.

The goal of this note is to propose a simple general framework for
justifying such convergence, applicable to both random measures and
random fields: a random field can conveniently be viewed as a random
vector in a Hilbert space with finite expectation of its squared norm; in this
case, the operators $Q^{(n_1,n_2)}$ have finite trace, and convergence
of their traces implies convergence of the random fields. Well-known
general arguments, see, for example, Theorem 4.2 in Billingsley's
book \cite{Billing}, justify interchanging limits in convergence in
distribution under the assumption of uniform convergence of the
approximants. This latter uniform convergence condition is particularly
convenient to verify in Sobolev spaces: indeed, a random field, for
example on the circle, defines a square-integrable random vector in an
appropriate Sobolev space, provided that the second moments of the
Fourier coefficients of our random field grow at most polynomially.
Thus, the trace-based framework for the kernels $Q^{(n_1,n_2)}$ makes
it possible to prove convergence in distribution of random fields.

\section{Random Elements of a Hilbert Space}
\subsection{The Sesquilinear Form Associated with a Pair of Random Elements of a Hilbert Space}
Let $H$ be a complex Hilbert space.
Fix an orthonormal basis $(e_n)_{n\in\NN}$ in $H$.

Now let $t_1,t_2$ be random vectors defined on a standard Borel
probability space $(\Omega,\mathcal{B},\PP)$ and taking values in $H$,
with
\[
  \EE_{\PP}\lvert t_1\rvert^{2}<+\infty,\qquad
  \EE_{\PP}\lvert t_2\rvert^{2}<+\infty,
\]
that is, our vectors $t_1,t_2$ are square-integrable.
Associate with the pair of random vectors $(t_1,t_2)$ the sesquilinear
form $Q^{(t_1,t_2)}$ on $H$ defined by
\[
  Q^{(t_1,t_2)}(h_1,h_2)
  =\EE_{\PP}\,\langle h_1,t_1\rangle\,\overline{\langle h_2,t_2\rangle}
  =\int_{\Omega}\langle h_1,t_1(\omega)\rangle\cdot\langle t_2(\omega),h_2\rangle\,d\PP(\omega).
\]
It follows directly from the definitions that
$$Q^{(t_1,t_2)}(h_1,h_2)=\overline{Q^{(t_2,t_1)}(h_2,h_1)},$$
while the form $Q^{(t,t)}$ is Hermitian and positive semidefinite.

A bounded sesquilinear form $Q$ on a Hilbert space is naturally
identified with a linear operator, which, by a slight abuse of notation,
we denote by the same symbol $Q$, according to the formula
\[
  \langle Qh_1,h_2\rangle=Q(h_1,h_2).
\]
The Hermitian property of the form is equivalent to the self-adjointness
of the operator, and positivity of the form is equivalent to positivity
of the operator.

% --- Here comes the block with (Q2)--(Q3) from the previous file ---
It follows directly from the definitions that

\begin{proposition} For any square-integrable random vectors
$t_1,t_2$ defined on a standard Borel probability space
$(\Omega,\mathcal{B},\PP)$, the corresponding operator
$Q^{(t_1,t_2)}$ is trace class, and we have
\[
  \operatorname{tr}Q^{(t_1,t_2)}=\EE_{\PP}\langle t_2,t_1\rangle.
\]
\end{proposition}

\begin{proof} First suppose that $t_1=t_2=t$.
Then the operator $Q^{(t,t)}$ is Hermitian and positive semidefinite,
and, for any orthonormal basis $(e_n)_{n\in\NN}$ in $H$, we have
\[
  \EE_{\PP}\langle t,t\rangle
  =\EE_{\PP}\sum_{n=1}^{+\infty}\bigl\lvert\langle t,e_n\rangle\bigr\rvert^{2}
  =\sum_{n=1}^{+\infty}\bigl\langle Q^{(t,t)}e_n,e_n\bigr\rangle
  =\operatorname{tr}Q^{(t,t)}.
\]
Applying this identity to the random vectors $t_1+t_2$ and
$t_1+it_2$ completes the proof of the proposition. \end{proof}

Next, by definition, we have
\[
  \EE_{\PP}\lvert t_1-t_2\rvert^{2}
  =\operatorname{tr}\Bigl(
     Q^{(t_1,t_1)}-Q^{(t_1,t_2)}-Q^{(t_2,t_1)}+Q^{(t_2,t_2)}
   \Bigr).
\]

% --- Mean-square convergence via traces ---
It is therefore clear that if $t^{(n)}$ is a sequence of random vectors
defined on a standard Borel probability space
$(\Omega,\mathcal{B},\PP)$ and taking values in a Hilbert space $H$,
with $\EE_{\PP}\lvert t^{(n)}\rvert^{2}<+\infty$ for every
$n\in\NN$, and if the sesquilinear forms
\[
  Q^{(kl)}(h_1,h_2)=\EE_{\PP}\,\langle h_1,t^{(k)}\rangle\,\langle t^{(l)},h_2\rangle
\]
satisfy the limit relation
\[
  \lim_{k,l\to\infty}\operatorname{tr}Q^{(kl)}=q_0,
\]
then the random vectors $t^{(n)}$ converge to a limiting
random vector $t^{(\infty)}$ in mean square, in the sense that
\[
  \lim\limits_{n\to\infty}\EE_{\PP}\lvert t^{(n)}-t^{(\infty)}\rvert^{2}=0,
\]
and, moreover, we have
\[
  \EE_{\PP}\lvert t^{(\infty)}\rvert^{2}=q_0.
\]

Now suppose that, for every $m\in\NN$, we are given a sequence
$t_m^{(n)}$ of random vectors taking values in $H$ and defined on a
standard Borel probability space $(\Omega_m,\mathcal{B}_m,\PP_m)$,
satisfying the relation
\[
  \EE\,\lvert t_m^{(n)}\rvert^{2}<+\infty,
\]
and converging in mean square to a limiting random vector
$t_m^{(\infty)}$.

% --- The blue page: four statements ---
Our next step is to show that if
\begin{itemize}
\item[1)] the convergence of $t_m^{(n)}$ to $t_m^{(\infty)}$ is uniform in $m$; and
\item[2)] for every $n\in\NN$, the vectors $t_m^{(n)}$ converge
in distribution to some limiting vector $t_\infty^{(n)}$,
\end{itemize}
then
\begin{itemize}
\item[1)] the vectors $t_m^{(\infty)}$ converge, as $m\to\infty$,
in distribution to some square-integrable random vector
$t_\infty^{(\infty)}$; and
\item[2)] the random vectors $t_\infty^{(n)}$ converge, as
$n\to+\infty$, in distribution to the vector $t_\infty^{(\infty)}$.
\end{itemize}

More formally, let $Q_m^{kl}$ stand for the sesquilinear form associated
with the vectors $t_m^{(k)}$ and $t_m^{(l)}$ and write:
\[
  Q_m^{kl}(h_1,h_2)
  =\EE_{\PP_m}\bigl\langle h_1,t_m^{(k)}\bigr\rangle
   \cdot\bigl\langle t_m^{(l)},h_2\bigr\rangle.
\]
We assume that
\[
  \sup_{m\in\NN,\,k\in\NN}\operatorname{tr}Q_m^{(k,k)}
  =\sup_{m\in\NN,\,k\in\NN}\EE_{\PP_m}\bigl\lvert t_m^{(k)}\bigr\rvert^{2}<+\infty.
\]
Our main assumption, namely the uniform convergence of the traces
$\operatorname{tr}Q_m^{(kl)}$, now takes the form:
\[
  \lim_{n\to+\infty}\sup_{m\in\NN}\sup_{k,l\ge n}
  \Bigl\lvert\operatorname{tr}\bigl(Q_m^{(k,k)}-Q_m^{(kl)}\bigr)\Bigr\rvert=0.
\]

Under our assumptions, it is clear that
\[
  \lim_{n\to\infty}\sup_{m\in\NN}\sup_{k,l\ge n}
  \EE_{\PP_m}\bigl\lvert t_m^{(k)}-t_m^{(l)}\bigr\rvert^{2}=0,
\]
whence, for every $m\in\NN$, there exists a random vector
$t_m^{(\infty)}$ defined on the probability space
$(\Omega_m,\mathcal{B}_m,\PP_m)$ and square-integrable in the sense that
$\EE_{\PP_m}\lvert t_m^{(\infty)}\rvert^{2}<+\infty$, such that
\[
  \lim_{n\to+\infty}\sup_{m\in\NN}\sup_{k\ge n}
  \EE_{\PP_m}\bigl\lvert t_m^{(\infty)}-t_m^{(k)}\bigr\rvert^{2}=0.
\]
It is also clear that
\[
  \sup_{m\in\NN}\EE_{\PP_m}\bigl\lvert t_m^{(\infty)}\rvert^{2}<+\infty.
\]

Let us make a simple observation. If $\xi_1,\xi_2$ are two random
elements defined on a common standard Borel probability space
$(\Omega,\mathcal{B},\PP)$ and taking values in a complete separable
metric space $S$, and if
\[
  d_2(\xi_1,\xi_2)
  =\EE_{\PP}\,\bigl\lvert d_S(\xi_1,\xi_2)\bigr\rvert^{2}
  =\int_{\Omega}\bigl\lvert d_S(\xi_1(\omega),\xi_2(\omega))\bigr\rvert^{2}\,d\PP(\omega)
  <+\infty,
\]
then the L\'evy--Prokhorov distance (cf. \cite{Zolot}, \cite{Bogachev})
between the distributions of the random elements $\xi_1,\xi_2$
satisfies the following inequality
\[
  \dLP\bigl((\xi_1)_{*}\PP,(\xi_2)_{*}\PP\bigr)
  \le\bigl(d_2(\xi_1,\xi_2)\bigr)^{1/3},
\]
which follows directly from the definition of the L\'evy--Prokhorov
distance. In particular, for random vectors $t_1,t_2$ defined on a
common probability space $(\Omega,\mathcal{B},\PP)$ and taking values
in a Hilbert space $H$, we have
\[
  \dLP\bigl((t_1)_{*}\PP,(t_2)_{*}\PP\bigr)
  \le\Bigl(\EE_{\PP}\,\lvert t_1-t_2\rvert^{2}\Bigr)^{1/3}.
\]

For $n\in\NN\cup\{\infty\}$, denote $\eta_m^{(n)}=(t_m^{(n)})_{*}\PP_m$. By definition, $\eta_m^{(n)}$ is a Borel
probability measure on $H$, and we have
\[
  \sup_{m\in\NN,\,n\in\NN}\int_{H}\lvert h\rvert^{2}\,d\eta_m^{(n)}(h)<+\infty.
\]
Denote this supremum by $\sigma_0$.
From the relation
\[
  \lim_{n\to\infty}\sup_{m\in\NN}\sup_{k\ge n}
  \EE_{\PP_m}\bigl\lvert t_m^{(\infty)}-t_m^{(k)}\bigr\rvert^{2}=0
\]
we obtain
\[
  \lim_{n\to\infty}\sup_{m\in\NN}\sup_{k\ge n}
  \dLP\bigl(\eta_m^{(k)},\eta_m^{(\infty)}\bigr)=0.
\]

Now suppose that, for every $n$, the measures $\eta_m^{(n)}$
converge weakly in $\cM_{\mathrm{prob}}(H)$ to a limiting measure
$\eta_\infty^{(n)}$. We will show that:
\begin{itemize}
\item[1)] the sequence of measures $\eta_\infty^{(n)}$ converges
weakly, as $n\to\infty$, to a limiting measure
$\eta_\infty^{(\infty)}$;
\item[2)] the sequence of measures $\eta_m^{(\infty)}$ converges
weakly, as $m\to\infty$, to the same limiting measure
$\eta_\infty^{(\infty)}$;
\end{itemize}

% --- Continuation of the proof: convergence of the measures eta_m^(infty) ---
Indeed, let $\varepsilon>0$. Choose $n_0$ in such a way that
\[
  d\bigl(\eta_m^{(k)},\eta_m^{(l)}\bigr)\le\varepsilon/3
\]
for all $m\in\NN$ and all $k,l\ge n_0$, including $k,l=\infty$.
Then we also have the inequality
\[
  d\bigl(\eta_\infty^{(k)},\eta_\infty^{(l)}\bigr)\le\varepsilon,
\]
and hence the measures $\eta_\infty^{(n)}$ do indeed converge to
a limiting measure $\eta_\infty^{(\infty)}$. Next, passing to the
limit as $l\to\infty$, we obtain, for all $k\ge n_0$, the inequality
\[
  \sup_{m\in\NN}\,d\bigl(\eta_m^{(k)},\eta_m^{(\infty)}\bigr)\le\varepsilon.
\]

Now, for every $k\in\NN$, the measures $\eta_m^{(k)}$ converge
as $m\to+\infty$ and, therefore, there exists
$m_0=m_0(k)\in\NN$ such that
\[
  \dLP\bigl(\eta_\infty^{(k)},\eta_m^{(k)}\bigr)\le\varepsilon
\]
for all $m>m_0$. It follows that
\[
  \dLP\bigl(\eta_\infty^{(k)},\eta_m^{(\infty)}\bigr)\le2\varepsilon
\]
for all $k> n_0$ and $m> m_0(k)$, which means precisely that the
sequence $\eta_m^{(\infty)}$ also converges to $\eta_\infty^{(\infty)}$
as $m\to+\infty$.

Finally, it is clear that
\[
  \int_{H}\lvert h\rvert^{2}\,d\eta_\infty^{(\infty)}(h)\le \sigma_0.
\]

% --- Conclusion and a new assumption ---
Summing up, we have obtained a certain strengthening of Theorem 4.2
in Billingsley's book \cite{Billing}, which is convenient for our
purposes. As before, suppose that, for every $m\in\NN$, a sequence of
square-integrable random vectors $t_m^{(n)}$ taking values in a Hilbert
space $H$ is given on a standard Borel probability space
$(\Omega_m,\mathcal{B}_m,\PP_m)$, together with the sesquilinear forms
$Q_m^{kl}$ associated with these vectors, as above, by the formula
\[
  Q_m^{kl}(h_1,h_2)
  =\EE_{\PP_m}\bigl\langle h_1,t_m^{(k)}\bigr\rangle
   \cdot\bigl\langle t_m^{(l)},h_2\bigr\rangle.
\]

\begin{proposition}\label{prop-prelim} Suppose that
\begin{itemize}
\item[1)] $\sup_{m\in\NN}\sup_{k\in\NN}\operatorname{tr}Q_m^{kk}<+\infty$;
\item[2)] $\lim_{n\to\infty}\sup_{m\in\NN}\sup_{k,l\ge n}
\Bigl\lvert\operatorname{tr}\bigl(Q_m^{kk}-Q_m^{kl}\bigr)\Bigr\rvert=0$;
\item[3)] for every $k\in\NN$, the sequence of vectors $t_m^{(k)}$
converges in distribution, as $m\to\infty$, to a limiting distribution
$\eta_\infty^{(k)}$ on $H$,
\end{itemize}
then
\begin{itemize}
\item[1)] the distributions $\eta_\infty^{(k)}$ converge weakly, as
$k\to+\infty$, to a limiting distribution $\eta_\infty^{(\infty)}$;
\item[2)] the random vectors $t_m^{(\infty)}$, defined for every
$m\in\NN$ by the condition
$\lim_{n\to\infty}\EE_{\PP_m}\lvert t_m^{(n)}-t_m^{(\infty)}\rvert^{2}=0$,
converge in distribution to a limiting random vector
$t_\infty^{(\infty)}$;
\item[3)] we have
$$\displaystyle\int\limits_{H}\lvert h\rvert^{2}\,d\eta_\infty^{(\infty)}(h)<+\infty.$$
\end{itemize} \end{proposition}

Proposition \ref{prop-prelim} is somewhat stronger than
Theorem 4.2 in Billingsley's book \cite{Billing}. Proposition \ref{prop-prelim} shows that,
in Sobolev spaces, the condition of uniform convergence of traces
reduces to checking convergence of finite-dimensional distributions
and polynomial growth of the Fourier coefficients.

\subsection{Measures Supported on Sobolev Subspaces}

The framework of this note is based on a simple observation:
convergence of traces is easy to verify in a Sobolev space.
We illustrate this approach with an elementary observation: if a
measure in a Hilbert space is supported on a ball in a Sobolev
subspace, then the distance from our measure to its finite-dimensional
projection is uniformly bounded and tends to zero uniformly as the
dimension of the finite-dimensional subspace increases.

To the basis $(e_n)_{n\in\NN}$, associate a scale of dense Sobolev
subspaces of $H$: for every $\alpha>0$, we consider the dense
subspace
\[
  H^{(\alpha)}=\Bigl\{h\in H:\ \sum_{n\in\NN}n^{2\alpha}\,
  \bigl|\langle h,e_n\rangle\bigr|^{2}<+\infty\Bigr\}.
\]
With respect to the Hilbert norm $\lVert\cdot\rVert_{H^{(\alpha)}}$
defined by the formula
\[
  \lVert h\rVert_{H^{(\alpha)}}^{2}=\sum_{n\in\NN}n^{2\alpha}\,
  \bigl|\langle h,e_n\rangle\bigr|^{2},
\]
the subspace $H^{(\alpha)}$ is itself a Hilbert space.

Next, let $\cM(H)$ be the space of Borel probability measures on $H$,
endowed with the L\'evy--Prokhorov distance $d_{LP}$, which induces
the weak topology. For $\alpha>0$ and $R>0$, denote by
$\cM_R^{(\alpha)}(H)$ the subset of measures $\PP\in\cM(H)$ such that
\[
  \EE_{\PP}\,\lVert h\rVert_{H^{(\alpha)}}^{2}\le R^2.
\]
Let $\Pi_N$ be the projection onto the finite-dimensional space spanned
by the vectors $(e_n)_{n\leq N}$, and let $\Pi_N^{\perp}=1-\Pi_N$ be
the projection onto the closed linear span of the vectors $(e_n)_{n> N}$.
From the definitions we directly obtain

%--- (2) ---
\begin{proposition}
For any $\alpha > 0$, $R > 0$, any measure
$\PP \in \cM_R^{(\alpha)}(H)$, and any $N \in \NN$, we have
\[
  \EE_{\PP}\, \lVert \Pi_N^{\perp} h \rVert_H^2 \le \frac{R^2}{N^{2\alpha}}.
\]\end{proposition}

\begin{corollary}
For any $\alpha > 0$, $R > 0$, any measures
$\PP_1, \PP_2 \in \cM_R^{(\alpha)}(H)$, and any $n \in \NN$, we have
\[
  d_{LP}(\PP_1,\PP_2) \le d_{LP}\bigl( (\Pi_n)_* \PP_1,\, (\Pi_n)_* \PP_2 \bigr)
  + \frac{2R^{2/3}}{n^{2\alpha/3}}.
\]\end{corollary}

\subsection{A Scale of Coordinate Sobolev Spaces}
Now consider a sequence of variables $(v_n)_{n \in \NN}$ and the
linear space $V(0)$ of finite sums:
\[
  V(0) = \Bigl\{ v = \sum_{k=1}^{N} \alpha_k v_k,\ \alpha_k \in \CC,\ N \in \NN \Bigr\}.
\]
We also associate with the sequence of variables $(v_n)_{n\in\NN}$
a Sobolev scale of Hilbert spaces: for every $\alpha\in \mathbb R$,
we set
\[
  HV^{\alpha} = \Bigl\{ v = \sum_{k=1}^{+\infty} \alpha_k v_k,\ \alpha_k \in \CC,\
  \sum_{k=1}^{+\infty} \bigl(k^{2\alpha}\bigr) |\alpha_k|^2 < +\infty \Bigr\}
\]
and write
$
  \| v \|_{HV^{\alpha}}^2 = \sum_{k=1}^{+\infty} \bigl( k^{2\alpha}\bigr) |\alpha_k|^2.
$

Now suppose that, to each variable $v_k$, a random
variable $\rX(v_k)$ is assigned on a standard Borel probability space
$(\Omega, \PP)$. To $v \in V(0)$, $v = \sum_{k=1}^{N} \alpha_k v_k$,
we then assign the random variable
$$
  \rX(v) = \sum_{k=1}^{N} \alpha_k\, \rX(v_k).
$$
We give simple sufficient conditions under which the correspondence
$v \to \rX(v)$ extends from $V(0)$ to $HV^{-\alpha}$.

%--- (5) ---
\begin{proposition} Let $\alpha\in\mathbb R$, $\delta_0 > 0$ be such that we have
\[
  \sup_{k \in \NN}
  k^{2\alpha+1+\delta_0}\EE\, |\rX(v_k)|^2 < +\infty.
\]
Then
\begin{enumerate}
  \item for $\PP$-almost every $\omega$, the random vector
    $
      t^{\rX} = \sum_{k=1}^{+\infty} \overline{ \rX(v_k)}\, v_k
    $
    belongs to the space $HV^{\alpha}$ and
    $$
    \EE_{\PP} \lVert t^{\rX} \rVert_{HV^{\alpha}}^2<+\infty;
    $$
  \item a random variable $\rX(v)$, defined by the formula
    $$
      \rX(v) = \langle v,\, t^{\rX} \rangle_{HV^{0}},
    $$
    may be assigned to every vector $v \in HV^{-\alpha}$, and we have
    $$
    \EE_{\PP} |\rX(v)|^2 \leq  \EE\, \lVert t^{\rX} \rVert_{HV^{\alpha}}^2 \lVert v \rVert_{HV^{-\alpha}}^2.
    $$
\end{enumerate} \end{proposition}
To the family $\bigl(\mathscr{X}(v_k)\bigr)_{k\in\NN}$ of random
variables associated with the variables $(v_k)_{k\in\NN}$, we
assign the sesquilinear form $Q^{\mathscr{X}}$ on the space $V(0)$,
defined by the formula
\[
  Q^{\mathscr{X}}\Bigl(\sum_{k=1}^{n}\alpha_k v_k,\,
  \sum_{l=1}^{m}\beta_l v_l\Bigr)
  =\sum_{k,l}\alpha_k\overline{\beta_l}\,
   \EE_{\PP}\,\mathscr{X}(v_k)\,\overline{\mathscr{X}(v_l)}.
\]
The form $Q^{\mathscr{X}}$ is Hermitian and positive semidefinite on
$V(0)$. Denote
\[
  Q^{\mathscr{X}}_{kl}
  =\EE_{\PP}\,\mathscr{X}(v_k)\,\overline{\mathscr{X}(v_l)}
  =Q^{\mathscr{X}}(v_k,v_l).
\]

It follows directly from the definitions that, if $\alpha\in\RR$ and
\[
  \sum_{k=1}^{+\infty} \frac{Q_{kk}}{k^{2\alpha}}<+\infty,
\]
then the form $Q^{\mathscr{X}}$ extends by continuity to a
sesquilinear Hermitian form with finite trace on $HV^{\alpha}$.
In particular, if the diagonal entries $Q_{kk}^{\mathscr{X}}$ have
polynomial growth, that is, if
\[
  \sup_{k\ge1}\frac{\log (2+Q^{\mathscr{X}}_{kk})}{\log(2+k)}<+\infty,
\]
then there exists $\alpha\in\RR$ such that the form
$Q^{\mathscr{X}}$ extends by continuity to a sesquilinear Hermitian
form with finite trace on $HV^{\alpha}$:
\begin{proposition}
Suppose that $Q_{kl}$ is a positive semidefinite Hermitian matrix
($Q_{kl}=\overline{Q_{lk}}$), and that
\[
  \sum_{k=1}^{+\infty}\frac{Q_{kk}}{k^{2\alpha}}=R<+\infty.
\]
Then the matrix $$Q_{kl}^{(\alpha)}=\frac{Q_{kl}}{k^{\alpha}l^{\alpha}}$$
is the matrix of a self-adjoint operator with trace $R$ on
$HV^{\alpha}$.
\end{proposition}
\begin{proof}\quad Set $$e_k^{(\alpha)}=\frac{v_k}{k^{\alpha}},\ k\in\NN.$$ 
We first prove the proposition for finitely supported matrices $Q$;
indeed, a finitely supported matrix $Q_{kl}$ defines a bounded
operator $Q^{(\alpha)}$ on $HV^{\alpha}$, given by the formula
\[
  Q^{(\alpha)}v_k=\sum_{l=1}^{+\infty}\frac{Q_{kl}}{l^{2\alpha}}\,v_l,
\]
or, equivalently, by the formula
\[
  Q^{(\alpha)}e_k^{(\alpha)}
  =\sum_{l=1}^{+\infty}\frac{Q_{kl}}{k^{\alpha}l^{\alpha}}\,e_l^{(\alpha)}.
\]
To the matrix $Q$, we assign a sesquilinear form on $V(0)$ defined by the formula
\[
  Q(\sum\limits_{k=1}^{N} c_k v_k,\sum\limits_{l=1}^{N} d_l v_l)=\sum_{k,l=1}^{N}Q_{kl}\,c_k\overline{d_l}.
\]
By definition, we have
\[
  \bigl\langle Q^{(\alpha)}e_k^{(\alpha)},e_l^{(\alpha)}\bigr\rangle_{HV^{\alpha}}
  =Q\bigl(e_k^{(\alpha)},e_l^{(\alpha)}\bigr);
\]
it is also clear that the operator $Q^{(\alpha)}$ is self-adjoint and
positive on $HV^{\alpha}$.
We have
\[
  \operatorname{tr}Q^{(\alpha)}
  =\sum_{k=1}^{+\infty}
   \bigl\langle Q^{(\alpha)}e_k^{(\alpha)},e_k^{(\alpha)}\bigr\rangle_{HV^{\alpha}}
  =\sum_{k=1}^{+\infty}\frac{Q_{kk}}{k^{2\alpha}}.
\]
Thus, the proposition is proved for finitely supported matrices.

Now, approximate the matrix $Q$ by positive semidefinite Hermitian matrices
$Q_N$ defined by
\[
  (Q_N)_{kl}=Q_{kl}\cdot\mathbb{I}_{\{k\le N,\,l\le N\}}.
\]
It is clear that, for $v,w\in V(0)$, the limit
$$\lim_{N\to\infty}\langle Q_N^{(\alpha)}v,w\rangle_{HV^{\alpha}}.$$
exists. It is also clear that
$$\lim_{N\to+\infty}\operatorname{tr}^{HV^{\alpha}}Q_N^{(\alpha)}=R.$$
It follows that the operators $Q_N^{(\alpha)}$ converge, as
$N\to+\infty$, in the class of trace-class operators acting on the
space $HV^{\alpha}$. It is also clear that the matrix of the limiting
operator is precisely $Q^{(\alpha)}$. The proposition is proved.
\end{proof}

\section{Random Fields on the Circle}

We now turn to Sobolev-regular random fields
on the unit circle.
Here we consider only the standard scale of Sobolev spaces.

Suppose that, for every $m \in \mathbb{N} \cup \{\infty\}$, a
sequence of random fields $\mathscr{X}_m^{(n)}$, $n \in \mathbb{N}$,
is given on a standard Borel probability space
$(\Omega_m, \mathcal{B}_m, \mathbb{P}_m)$, and that, for every
$k \in \mathbb{Z}$, we have
\[
\mathscr{X}_m^{(n)}(\exp(ik\theta)) \in L_2(\Omega_m,\mathbb{P}_m), \quad m, n \in \mathbb{N} \cup \{\infty\}.
\]

As before, consider the correlation matrices
$Q_m^{(n_1, n_2)}(k_1, k_2)$, defined for
$m, n_1, n_2 \in \mathbb{N} \cup \{\infty\}$ and
$k_1, k_2 \in \mathbb{Z}$ by the formula
\[
Q_m^{(n_1, n_2)}(k_1, k_2) = \mathbb{E}_{\mathbb{P}_m}\,
\mathscr{X}_m^{(n_1)}(\exp(ik_1\theta)) \cdot
\overline{\mathscr{X}_m^{(n_2)}(\exp(ik_2\theta))}.
\]

The following assumption of uniform polynomial boundedness of the
matrix coefficients $Q_m^{(n,n)}(k,k)$ plays a key role for us.

\begin{assumption}\label{poly-bnd}
There exist $C > 0$ and $\alpha > 0$ such that, for all
$m, n \in \mathbb{N} \cup \{\infty\}$ and $k \in \mathbb{Z}$, the
following inequality holds:
\[
Q_m^{(n,n)}(k,k) \leq C\,(1 + |k|)^\alpha.
\]
\end{assumption}

Under Assumption \ref{poly-bnd}, every random field
$\mathscr{X}_m^{(n)}$, $m, n \in \mathbb{N} \cup \{\infty\}$, induces a square-integrable random vector $t_m^{(n)}$ defined
on the probability space $(\Omega_m, \mathcal{B}_m, \mathbb{P}_m)$, and which is square-integrable in the Sobolev space $H^{(-\alpha-1)}(\mathbb{T})$:
for all $m, n \in \mathbb{N} \cup \{\infty\}$, we have
\[
\mathbb{E}\,\bigl\|t_m^{(n)}\bigr\|^2_{H^{(-\alpha-1)}(\mathbb{T})} < +\infty
\]
and, moreover,
\[
\sup_{m, n \in \mathbb{N} \cup \{\infty\}} \mathbb{E}\,\bigl\|t_m^{(n)}\bigr\|^2_{H^{(-\alpha-1)}(\mathbb{T})} < +\infty.
\]

Proposition \ref{prop-prelim} now yields the following

\begin{corollary}\label{cor-poly}
Suppose that, for every $r \in \mathbb{N}$ and every
$n_1, \ldots, n_r \in \mathbb{N}$, $k_1, \ldots, k_r \in \mathbb{Z}$,
the sequences of random vectors
\[
\bigl(\mathscr{X}_m^{(n_1)}(\exp(ik_1\theta)), \ldots,
\mathscr{X}_m^{(n_r)}(\exp(ik_r\theta))\bigr)
\]
converge in distribution, as $m \to \infty$, to the random vector
\[
\bigl(\mathscr{X}_\infty^{(n_1)}(\exp(ik_1\theta)), \ldots,
\mathscr{X}_\infty^{(n_r)}(\exp(ik_r\theta))\bigr).
\]

Furthermore, suppose that, for every
$m \in \mathbb{N} \cup \{\infty\}$ and every
$k_1, k_2 \in \mathbb{Z}$, the joint limit
\[
\lim_{n_1, n_2 \to \infty} Q_m^{(n_1, n_2)}(k_1, k_2)
= Q_m^{(\infty, \infty)}(k_1, k_2),
\]
exists, and that
$$
\lim_{n\to\infty} \sup\limits_{m\in \mathbb {N}\cup \{\infty\}}
\sup\limits_{n_1, n_2 \in [n, +\infty]}
|Q_m^{(n_1, n_2)}(k_1, k_2) - Q_m^{(\infty, \infty)}(k_1, k_2)|=0.
$$

Then:
\begin{enumerate}
    \item for every $m \in \mathbb{N} \cup \{\infty\}$, the sequence
    of random vectors $t_m^{(n)}$, taking values in
    $H^{(-\alpha-1)}(\mathbb{T})$, converges in
    $H^{(-\alpha-1)}(\mathbb{T})$ to a limiting random vector
    $t_m^{(\infty)}$:
    \[
    \mathbb{E}_{\mathbb{P}_m}\,\bigl\|t_m^{(n)} - t_m^{(\infty)}\bigr\|^2_{H^{(-\alpha-1)}(\mathbb{T})} \to 0;
    \]
    \item the vectors $t_m^{(\infty)}$ converge
    in distribution to the vector $t_\infty^{(\infty)}$ as $m \to \infty$.
\end{enumerate}
\end{corollary}

\section{Double Averaging Lemma}

\subsection{Diagonal Approximations}
Let $(W, d_W)$ be a compact metric space endowed with a nonatomic
Borel probability measure $\eta$.

\begin{definition}
A sequence of nonnegative bounded Borel functions
\[
P_n : W \times W \to [0, +\infty)
\]
will be called a \textit{diagonal approximation} if:
\begin{enumerate}
    \item for every $s \in W$ and $n \in \mathbb{N}$, we have
    \[
    \int_W P_n(s, t)\, d\eta(t) = 1;
    \]
    \item the following inequality holds:
    \[
    \sup_{n \in \mathbb{N}} \sup_{t \in W} \int_W P_n(s, t)\, d\eta(s) < +\infty;
    \]
    \item there exists a positive sequence $\varepsilon_n > 0$,
    $\varepsilon_n \to 0$, such that $P_n(s, t) = 0$ whenever
    $d_W(s, t) > \varepsilon_n$.
\end{enumerate}
\end{definition}

\begin{definition}
A Borel function $Q$ on $W \times W$ will be called
\textit{uniformly tempered along the diagonal} if:
\begin{enumerate}
    \item the restriction of $Q$ to the subset
    \[
    W \times W \setminus \mathrm{diag} =
    \{(s, t) \in W \times W : s \neq t\};
    \]
    is continuous;
    \item the following relation holds:
    \[
    \lim_{\varepsilon \to 0} \sup_{s \in W}
    \int_{\{t : d_W(s, t) \leq \varepsilon\}} |Q(s, t)|\, d\eta(t) = 0.
    \]
\end{enumerate}
\end{definition}

Our assumptions imply the following simple result.

\begin{proposition}\label{temp-conv}
If $Q$ is uniformly tempered along the diagonal and $P_n$ is a
diagonal approximation, then, as $n \to \infty$, the sequence of
functions $Z_n$ defined by the formula
\[
Z_n(s_1, s_2) = \int\limits_W \int\limits_W
P_n(s_1, t_1)\, P_n(s_2, t_2)\, Q(t_1, t_2)\,
d\eta(t_1)\, d\eta(t_2),
\]
converges to $Q$ pointwise off the diagonal, in
$L_1(W \times W, \eta \otimes \eta)$, and uniformly outside any
$\varepsilon$-neighborhood of the diagonal.
\end{proposition}

Jensen's inequality allows us to use the formalism of diagonal approximations in our work with exponentials of uniformly tempered
functions.

Let $K$ be a Borel function on $W \times W$, uniformly tempered along
the diagonal, and such that $\exp K$ is also uniformly tempered along
the diagonal. Let $P_n(s,t)$ be a diagonal approximation. Define the
sequence
\[
Y_n(s_1, s_2) = \exp \iint_{W \times W}
P_n(s_1, t_1)\, P_n(s_2, t_2)\, K(t_1, t_2)\,
d\eta(t_1)\, d\eta(t_2).
\]
Jensen's inequality immediately implies the following
\begin{proposition}\label{temp-conv-0}
As $n \to \infty$, the sequence $Y_n$ converges to $\exp(K)$
pointwise off the diagonal, in $L_1(W \times W, \eta \otimes \eta)$,
and uniformly outside any $\varepsilon$-neighborhood of the diagonal.
\end{proposition}

We now give a uniform version of Proposition \ref{temp-conv-0}.
Let $P^{(1)}_n(s,t)$ and $P^{(2)}_n(s,t)$ be diagonal approximations.
Define a bounded Borel function $R_{jl}(s,t)$ by
$$
R_{jl}(s,t)=\int\limits_{W} P^{(1)}_j(s, u)\, P_l^{(2)}(u, t)d\eta(u).
$$
As above, let $K$ be a Borel function on $W \times W$, uniformly
tempered along the diagonal, and such that $\exp K$ is also uniformly
tempered along the diagonal.
Denote
\[
Z_{nl}^{(1)}(s_1, s_2) = \int\limits_W \int\limits_W
P_n^{(1)}(s_1, t_1)\, P_l^{(1)}(s_2, t_2)\, K(t_1, t_2)\,
d\eta(t_1)\, d\eta(t_2),
\]
\[
Z_{nl}^{(2)}(s_1, s_2) = \int\limits_W \int\limits_W
P_n^{(2)}(s_1, t_1)\, P_l^{(2)}(s_2, t_2)\, K(t_1, t_2)\,
d\eta(t_1)\, d\eta(t_2),
\]
and introduce the sequence
\[
Y_{abcd}(s_1, s_2) = \exp \int\limits_{W}\int\limits_{W}
R_{ab}(s_1, t_1)\, R_{cd}(s_2, t_2)\, K(t_1, t_2)\,
d\eta(t_1)\, d\eta(t_2).
\]

Directly from the definitions, we obtain that
\begin{proposition}\label{temp-conv-11}\begin{enumerate}
\item As $a,c \to \infty$, the sequence $Y_{abcd}$ converges to
$\exp(Z_{bd}^{(2)})$ in $L_1(W \times W, \eta \otimes \eta)$,
uniformly in $b,d$;
\item As $b,d \to \infty$, the sequence $Y_{abcd}$ converges to
$\exp(Z_{ac}^{(1)})$ in $L_1(W \times W, \eta \otimes \eta)$,
uniformly in $a,c$;
\item If $q_{abcd}$ is a uniformly bounded family of Borel functions
on $W \times W$ which, uniformly in $b,d$, converges to $1$ on
$W \times W$ outside every neighborhood of the diagonal as
$a,c \to \infty$, then the sequence $q_{abcd}Y_{abcd}$ also converges,
as $a,c \to \infty$, to $\exp(Z_{bd}^{(2)})$ in
$L_1(W \times W, \eta \otimes \eta)$, uniformly in $b,d$. If,
on the other hand, $q_{abcd}$ converges to $1$ on $W \times W$ outside
every neighbourhood of the diagonal, uniformly in $a,c$, as
$b,d \to \infty$, then the sequence $q_{abcd}Y_{abcd}$ converges to
$\exp(Z_{ac}^{(1)})$ in $L_1(W \times W, \eta \otimes \eta)$,
uniformly in $a,c$.
\end{enumerate}
\end{proposition}

Proposition \ref{temp-conv-11} directly implies the main result of this section, the following Double-Averaging Lemma.

\begin{lemma}\label{twav}
\begin{enumerate}
\item As $a,b,c,d \to \infty$, the sequence $Y_{abcd}$ converges, in the totality of its indices, to $\exp(K)$ pointwise off
the diagonal, in $L_1(W \times W, \eta \otimes \eta)$, and uniformly
outside any $\varepsilon$-neighbourhood of the diagonal.
\item If $q_{abcd}$ is a uniformly bounded family of Borel functions
on $W \times W$ which converges uniformly, in the totality of its indices, to $1$ on $W \times W$ outside every neighbourhood of the
diagonal, then the sequence $q_{abcd}Y_{abcd}$ also converges, in the totality of its indices, to $\exp(K)$ pointwise off the
diagonal, in $L_1(W \times W, \eta \otimes \eta)$, and uniformly
outside any $\varepsilon$-neighbourhood of the diagonal.
\end{enumerate}
\end{lemma}

\section{Convergence of Random Measures to Gaussian Multiplicative Chaos}

We apply the general framework developed above to prove convergence
of exponentials of random fields to Gaussian Multiplicative Chaos
(cf. \cite{Kahane}, \cite{Berestyczki}, \cite{Bufetov-UMN},
\cite{FyodorovBouchaud}, \cite{FyodorovHorSimm},\cite{Lacoin},
\cite{RhodesVargas}, \cite{Shamov}, \cite{Webb}). There are several
well-known results in the literature on the convergence of
exponentials of non-Gaussian fields to Gaussian Multiplicative Chaos;
see, for example, Berestycki, Webb, and Wong
\cite{bewewo}, Lambert, Ostrovsky, and Simm
\cite{lambert-ostrovsky-simm}, and Claeys, Fahs, Lambert, and Webb
\cite{claeys}. In the $L_2$ case, Proposition \ref{chaos-conv}
allows for a more general approximation scheme.
Here we restrict ourselves to the case of the circle and real-valued
chaos. Sufficient conditions for convergence to complex Gaussian
Multiplicative Chaos can also be obtained by the trace method and will be considered in the sequel to this paper.

As before, we write $\T=\RR/2\pi\ZZ$. A function $f\in L_1(\T)$ has Fourier coefficients
\[
  \widehat{f}(k)=\frac{1}{2\pi}\int_{\T}e^{-ik\theta}f(\theta)\,d\theta.
\]
If $f\in L_2(\T)$, then, of course, we have
\[
  f=\sum_{k\in\ZZ}\widehat{f}(k)e^{ik\theta},
\]
and the series on the right-hand side converges in $L_2(\T)$.

To the standard trigonometric basis, we assign the corresponding
scale of Sobolev spaces: for $\alpha\in \mathbb R$, define the norm
$\lVert f\rVert_{H^{\alpha}}$ by the formula
$$
\lVert f\rVert_{H^{\alpha}}^{2}=\sum\limits_{n\in \mathbb Z}
(1+n^2)^{\alpha} |\widehat{f}(n)|^2.
$$
% --- end of insertion ---

We start with a simple observation. Let $(\Omega, \mathcal{B}, \mathbb{P})$
be a standard probability space. Let
$\xi_1, \xi_2: \Omega \to H^{-1}(\mathbb{T})$ be square-integrable
random vectors. We will say that the pairwise correlations of the
random vectors $\xi_1, \xi_2$ are represented by the correlation
matrix $Q_{rj}$, $r,j\in\{1,2\}$, where
$Q_{rj} \in L_1(\mathbb{T} \times \mathbb{T})$, if, for all
$r,j\in\{1,2\}$ and $k, \ell \in \mathbb{Z}$, the following identity
holds:
\[
\mathbb{E}_{\mathbb{P}}\xi_r\bigl(\exp(ik\theta)\bigr)
\overline{\xi_j\bigl(\exp(i\ell\theta)\bigr)} =
\widehat{Q_{rj}}(-k, \ell).
\]
We now immediately have
\begin{proposition}\label{hone-lone}
\[
\mathbb{E}_{\mathbb{P}}\, \|\xi_1 - \xi_2\|^2_{H^{-1}(\mathbb{T})}
\leq 4\Bigl(\|Q_{11} - Q_{12}\|_{L_1(\mathbb{T}\times\mathbb{T})}
+ \|Q_{21} - Q_{22}\|_{L_1(\mathbb{T}\times\mathbb{T})}\Bigr).
\]
\end{proposition}
\begin{proof}
Indeed, we have:
\begin{multline}\label{qq}
\mathbb{E}_{\mathbb{P}}\, \|\xi_1 - \xi_2\|^2_{H^{-1}(\mathbb{T})} =\\= \sum_{n \in \mathbb{Z}} \frac{\mathbb{E}_{\mathbb{P}}\, |\widehat{\xi}_1(n) - \widehat{\xi}_2(n)|^2}{1 + n^2} =\\= \sum_{n \in \mathbb{Z}} \frac{\widehat{Q}_{11}(-n,n) - \widehat{Q}_{12}(-n,n) - \widehat{Q}_{21}(-n,n) + \widehat{Q}_{22}(-n,n)}{1 + n^2}\leq \\
\leq 4\Bigl(\|Q_{11} - Q_{12}\|_{L_1(\mathbb{T}\times\mathbb{T})} + \|Q_{21} - Q_{22}\|_{L_1(\mathbb{T}\times\mathbb{T})}\Bigr),
\end{multline}
as announced, and the Proposition \ref{hone-lone} is proved.
\end{proof}

Next, the Sobolev space $\dH$ is defined as the completion of the linear
space of trigonometric polynomials with zero constant term with respect
to the Hilbert norm $\lVert f\rVert_{\dH}$ given by the formula
\[
  \lVert f\rVert_{\dH}^{2}
  =\sum_{k\in\ZZ\setminus\{0\}}\frac{\lvert\widehat{f}(k)\rvert^{2}}{2\lvert k\rvert}
  =\int\limits_{\T} \int\limits_{\T}f(\theta_1)\overline{f(\theta_2)}\,
    \log\frac{1}{\lvert e^{i\theta_1}-e^{i\theta_2}\rvert}\,
    \frac{d\theta_1}{2\pi}\frac{d\theta_2}{2\pi}.
\]
Within the complex Sobolev space $\dH$, consider the real subspace
$\ddH$ spanned by real-valued trigonometric polynomials.
To the real Hilbert space $\ddH$ assign an isometric real
Gaussian Hilbert space $\cH$: in other words, for every
$k\in\mathbb N$, let $\xi_k^+, \xi_k^-$ be real Gaussian random
variables with mean $0$ and variance $1/4k$, in such a way that the
random variables $(\xi_k^{\pm})_{k\in\mathbb N}$ are jointly
independent. Denote by $\cH$ the closed real linear span of the
vectors $\xi^{\pm}_k$.

The map $$\cX_\infty\colon \cos(k\theta)\mapsto\xi_k^+,
\sin (k\theta)\mapsto \xi_k^-$$ defines an isometric isomorphism
between the Hilbert space $\dH^{\mathbb R}$ and the Gaussian Hilbert
space $\cH$ (cf. e.g. \cite{lifshits}); moreover, any
function $f\in\dH^{\mathbb R}$ with zero average, represented by its real Fourier series
$$
f=\sum_{k\in\mathbb N} a_k \cos(k\theta)+b_k \sin(k\theta),
$$
induces the Gaussian vector
\[
  \cX_\infty(f)=\sum_{k\in\mathbb N}a_k\xi_k^+ + b_k \xi_k^-,
\]
and we have $$\EE\lvert\cX_\infty(f)\rvert^{2}=\lVert f\rVert_{\dH}^{2}.$$ Finally, set $\cX_\infty(1)=0$.

Let us recall one possible definition of Gaussian Multiplicative
Chaos. In this note, we are concerned only with the $L_2$ phase; the
$L_1$ phase will be considered in the sequel to this paper.
Fix $\gamma \in (0,1)$. Here it is convenient for us to consider a
specific sequence of approximations to Gaussian Multiplicative Chaos,
setting aside the question of the independence of the limit from the
particular choice of approximation. Thus, fix a sequence of
partitions $\pi_n$ of our circle $\T$ into arcs of the form
$$\left[\frac{k}{2^n}\pi,\frac{k+1}{2^n}\pi\right),
\quad k=0,\dots,2^{n+1}-1,$$
and consider the sequence of stochastic processes
\[
  \cX_\infty^n(e^{i\theta})=\sum_{I\in\pi_n}\frac{\cX_\infty(\mathbf{1}_I)}{|I|}\,\ind_I(\theta)
\]
(where $|I|=\int_I\frac{d\theta}{2\pi}$ stands for the length of the arc $I$).
Consider, next, the sequence of normalized exponentials
\[
  dM^n_{GMC_\gamma}(\theta)
  =\frac{\exp(\gamma\cX_\infty^n(e^{i\theta}))}{\EE\exp(\gamma\cX_\infty^n(e^{i\theta}))}\,
   \frac{d\theta}{2\pi}
\]
and the sequence of corresponding kernels
\[
  Q^{n_1,n_2}_{GMC_\gamma}(e^{i\theta_1}, e^{i\theta_2})
  =\frac{\EE\exp(\gamma\cX_\infty^{n_1}(e^{i\theta_1})+\gamma\cX_\infty^{n_2}(e^{i\theta_2}))}
        {\EE\exp(\gamma\cX_\infty^{n_1}(e^{i\theta_1}))\cdot \EE\exp(\gamma\cX_\infty^{n_2}(e^{i\theta_2}))}.
\]

It follows from Propositions \ref{temp-conv-0} and \ref{hone-lone}
that the operators $Q^{n_1,n_2}_{GMC_\gamma}$ are trace class in the
Sobolev space $H^{-1}(\T)$, and that the limit
\[
  \lim_{n_1,n_2\to\infty}\operatorname{tr}Q^{n_1,n_2}_{GMC_\gamma}.
\]
exists.
Thus, the sequence of random measures
$dM^n_{GMC_\gamma}(\theta)$ converges in the Sobolev space
$H^{-1}(\T)$ to a limiting random measure, precisely the Gaussian
Multiplicative Chaos $GMC_\gamma$.

We now give sufficient conditions for the convergence of a sequence of
exponentials of random processes to Gaussian Multiplicative Chaos.

Thus, let $\cX_\infty$ be a Gaussian random field on $\mathbb{T}$ with
correlation function
\[
Q(e^{i\theta_1}, e^{i\theta_2}) = \log \frac{1}{|e^{i\theta_1} - e^{i\theta_2}|}.
\]

Next, suppose that, for every $m\in\NN$, a random process $X_m$ with
continuous sample paths on $\mathbb{T}$ is given on a standard Borel
probability space $(\Omega_m,\mathcal{B}_m,\PP_m)$.

To the stochastic process $X_m$ we assign a random measure $M_m$,
defined on the probability space $(\Omega_m,\mathcal{B}_m,\PP_m)$,
by the formula
\[
  dM_m(e^{i\theta})=
  \frac{\exp(\gamma X_m(e^{i\theta}))}{\EE\exp(\gamma X_m(e^{i\theta}))}\,
  \frac{d\theta}{2\pi}.
\]

If $\varphi$ is a bounded Borel function on the circle, then we write
$$
X_m(\varphi)=\frac 1{2\pi}\int_{\mathbb{T}}
X_m(\exp(i\theta))\varphi(\theta) d\theta.
$$
For simplicity, we also assume that
$$
\widehat X_m(0)=\frac 1{2\pi} \int_{\mathbb{T}}
X_m(\exp(i\theta))d\theta=0
$$
(this assumption can be substantially weakened).

We need an assumption of uniform boundedness of the exponential
moments of the processes $X_m$: we assume that there exists
$\gamma_1>\gamma$ such that, for every bounded real-valued Borel
function $\varphi$ with zero mean, the following estimate holds:
\begin{equation}\label{unif-bdd}
\sup\limits_{m\in {\mathbb N}\cup\{\infty\}}\mathbb{E}
\exp\!\left(\gamma_1  X_m(\varphi)\right)
\leq \exp\!\left(\frac{\gamma_1^2}{2} ||\varphi||^2_{\dH}\right).
\end{equation}

Furthermore, suppose that, for any bounded real-valued Borel functions
$\varphi_1, \dots, \varphi_r$, the following convergence in
distribution of finite-dimensional random vectors holds:
\[
\bigl({X}_m(\varphi_1), \ldots, {X}_m(\varphi_r)\bigr)\to
\bigl(\cX_\infty(\varphi_1), \ldots,
\cX_\infty(\varphi_r)\bigr).
\]

Next, suppose that we are given a diagonal approximation $P_m(s,t)$
on $\mathbb{T}\times \mathbb{T}$.
Finally, let $\mathcal{M}_1(\mathbb{T})$ be the space of Borel
probability measures on the circle with the weak topology, and let
$\mathcal{M}_3(\mathbb{T})$ be the space of positive Borel measures on
the circle with total mass at most three, also endowed with the weak
topology. Let $R_m(\eta): \mathcal{M}_3(\mathbb{T}) \to \mathbb{R}$ be
continuous functions which are uniformly bounded, and such that
$$
\lim\limits_{m\to\infty} \sup\limits_{\eta_1, \eta_2\in
\mathcal{M}_1(\mathbb{T}), \mathrm{supp}(\eta_1)\cap
\mathrm{supp}(\eta_2)=\emptyset}
|R_m(\eta_1+\eta_2)-R_m(\eta_1)-R_m(\eta_2)|=0.
$$
We use $dt$ to denote the normalized Lebesgue measure on the circle.
Our key assumption is that, for every finite positive measure $\eta$
on $\mathbb{T}$ with total mass at most three, the following relation holds:
\begin{multline}\label{key-assump}
\mathbb{E} \exp\!\left(\gamma \int_{\mathbb{T}} X_m(s)\, d\eta(s)\right) =\\= \exp\!\Big(\frac{\gamma^2}{2}
\int\limits_{\mathbb{T}}\int\limits_{\mathbb{T}}\int\limits_{\mathbb{T}}\int\limits_{\mathbb{T}}
P_m(s_1, t_1) \cdot P_m(s_2, t_2)
\\\log \frac{1}{|\exp(it_1) - \exp(i t_2)|}\,
d\eta(s_1)\, d\eta(s_2){dt_1}{dt_2} + R_m(\eta)\Big).
\end{multline}

Conditions of this kind can conveniently be verified for determinantal
point processes (cf. \cite{Bufetov-FAA}, \cite{Bufetov-UMN},
\cite{Gorbunov}, \cite{Bufetov-Cond}, \cite{Bufetov12},
\cite{Bufetov-Symm}).

As usual, by the distribution of a random measure defined on a
standard Borel probability space, we mean its pushforward on
the space of probability measures; in our particular case, on the
space of probability measures on the space of finite measures on the
circle. Convergence in distribution is, as usual, understood with
respect to the weak topology on the space of probability measures.

Under Assumption \ref{key-assump}, we have the following result.

\begin{proposition}\label{chaos-conv} As $m\to+\infty$, the random
measures $M_m$ converge in distribution to the Gaussian Multiplicative
Chaos $\mathrm{GMC}_\gamma$. \end{proposition}

\begin{proof}

The symbol $\pi_n$, as above, will denote the partition of our circle
$\T$ into half-open intervals
$$\bigl[\frac{k}{2^n}\pi,\frac{k+1}{2^n}\pi\bigr),\
k=0,\dots,2^{n+1}-1.$$ Associate with the random process $X_m$ the
following sequence of approximations
%
%\begin{multline}

\[
X_m^{(n)} = \sum_{I \in \pi_n}
X_m\!\left(\frac{\mathbf{1}_I}{|I|}\right) \cdot 1_I.
\]

%\end{multline}
In other words, $X_m^{(n)}$ is the conditional expectation of $X_m$
with respect to the $\sigma$-algebra generated by the partition
$\pi_n$.

Next, consider the sequence of random measures
\[
  dM_m^{(n)}(e^{i\theta})
  =\frac{\exp(\gamma X_m^{(n)}(e^{i\theta}))}
        {\EE_{\PP_m}\exp(\gamma X_m^{(n)}(e^{i\theta}))}\,
        \frac{d\theta}{2\pi}.
\]
To the measures $M_m^{(n_1)}$, $M_m^{(n_2)}$ we assign the kernel
\[
  Q^{n_1,n_2}_{(m)}(e^{i\theta_1},e^{i\theta_2})
  =\frac{\EE_{\PP_m}\exp\bigl(\gamma (X_m^{(n_1)}(e^{i\theta_1}))
        +\gamma X_m^{(n_2)}(e^{i\theta_2}))\bigr)}
        {\EE_{\PP_m}\exp(\gamma X_m^{(n_1)}(e^{i\theta_1}))\,
         \EE_{\PP_m}\exp(\gamma X_m^{(n_2)}(e^{i\theta_2}))}.
\]

The finite-dimensional distributions of the Fourier coefficients of
our random measures converge as a consequence of the convergence of
the finite-dimensional distributions of the fields themselves, the
assumption of uniform boundedness of exponential moments, and the
following simple observation: if $\xi_n$ are random vectors in
$\mathbb{R}^m$ such that $\xi_n \to \xi_\infty$ in distribution, and
$\Phi: \mathbb{R}^m \to \mathbb{R}^k$ is a continuous function, then the convergence
$\Phi(\xi_n) \to \Phi(\xi_\infty)$ holds in distribution as $n \to \infty$.
Moreover, if $\alpha>0$ and
\[
\sup_{n \in \mathbb{N} \cup \{\infty\}}
\mathbb{E}\, |\Phi(\xi_n)|^{1+\alpha} < +\infty,
\]
then the expectations also converge:
$$
\EE \Phi(\xi_n) \to \EE \Phi(\xi_\infty).
$$

Since $\gamma \in (0,1)$, our assumptions and Proposition
\ref{temp-conv} directly imply that each of the random measures
$M_m^{(n)}$ defines a square-integrable random vector $t_m^{(n)}$
in the Sobolev space $H^{-1}(\mathbb{T})$.

Denote by $t_m$ the random vector associated
    with the measure $M_m$ itself. It is also clear from Proposition \ref{temp-conv}, Proposition \ref{hone-lone}, and Lemma \ref{twav} that
    \[
    \lim_{n \to \infty} \mathbb{E}_{\mathbb{P}_m}\,
    \bigl\|t_m^{(n)} - t_m\bigr\|^2_{H^{-1}(\mathbb{T})} = 0,
    \]
    where the convergence is uniform in $m$.

It is also clear that, for every $r \in \mathbb{N}$, as
\(m\to\infty\), for fixed \(n\), the random vector
    \[
    \left(\int e^{ik\theta}\, dM_m^{(n)}(\theta)\right)_{|k| \leq r}
    \]
    converges in distribution to the random vector
    \[
    \left(\int e^{ik\theta}\, dM_{\mathrm{GMC}_\gamma}^{(n)}(\theta)\right)_{|k| \leq r}.
    \]

The Proposition \ref{chaos-conv} now follows from Corollary \ref{cor-poly}.

\end{proof}

\noindent\textbf{ Acknowledgements.} I am deeply grateful to
A.S. Arakcheev, M.A. Borovikov, D.V. Gnatyuk, S.M. Gorbunov, and
I.I. Zavolokin for useful discussions.

It is my pleasant duty to thank the referees whose detailed
comments have greatly helped me improve the exposition.

The author is a winner of the <<Leader>> competition of the <<BASIS>>
Foundation for the Development of Theoretical Physics and Mathematics,
and thanks the jury and sponsors.

This paper was partially written as part of the work of a research
group at the Sirius Mathematics Centre, and I am deeply grateful to
the Sirius Mathematics Centre for its warm hospitality.

\end{document}